\documentclass{amsart}

\newtheorem{theorem}{Theorem}[section]
\newtheorem{lemma}[theorem]{Lemma}
\newtheorem{conjecture}[theorem]{Conjecture}

\theoremstyle{definition}
\newtheorem{definition}[theorem]{Definition}

\theoremstyle{remark}

\numberwithin{equation}{section}

\begin{document}

\title{Symmetric Word Equations in Two Positive Definite
Letters}
\thanks{This research was conducted, in part, during the summer of 1999
at the College of William and Mary's Research Experiences for
Undergraduates program and was supported
by NSF REU grant DMS-96-19577.}%

\author{Christopher J. Hillar}
\address{Department of Mathematics, University of California, Berkeley, CA 94720.}
\email{chillar@math.berkeley.edu}
\thanks{The work of the first author is supported under a National
Science Foundation Graduate Research Fellowship.}

\author{Charles R. Johnson}
\address{Department of
Mathematics, College of William and Mary, Williamsburg, VA
23187-8795.} \email{crjohnso@math.wm.edu}



\begin{abstract}
A generalized word in two positive definite matrices $A$ and $B$
is a finite product of nonzero real powers of $A$ and $B$.
Symmetric words in positive definite $A$ and $B$ are positive
definite, and so for fixed $B$, we can view a symmetric word,
$S(A,B)$, as a map from the set of positive definite matrices into
itself. Given positive definite $P$, $B$, and a symmetric word,
$S(A,B)$, with positive powers of $A$, we define a symmetric word
equation as an equation of the form $S(A,B) = P$. Such an equation
is solvable if there is always a positive definite solution $A$
for any given $B$ and $P$. We prove that all symmetric word
equations are solvable. Applications of this fact, methods for
solution, questions about unique solvability (injectivity), and
generalizations are also discussed.
\end{abstract}

\maketitle

\section{Introduction}

A {\it generalized word} ({\it g-word}, for short) $W = W(A,B)$ in
two letters $A$ and $B$ is an expression of the form $$ W = A^{p_1
} B^{q_1 } A^{p_2 } B^{q_2 } \cdots A^{p_k } B^{q_k } A^{p_{k + 1}
} $$ in which the exponents $p_{i}$ and $q_{i}$ are real numbers
such that $p_{i},q_{i} \neq 0$, $i = 1,\ldots,k$, and $p_{k+1}$ is
an arbitrary real number. We call $k$ the {\it class number} of
$W$. The {\it reversal} of the g-word $W$ is $W^{*} = A^{p_{k + 1}
} B^{q_k } A^{p_k } \cdots B^{q_2 } A^{p_2 } B^{q_1 } A^{p_1 }$
and a g-word is {\it symmetric} if it is identical to its reversal
(in other contexts, the name "palindromic" is also used).  We will
call a g-word, $W$, $A$-{\it positive} ($A$-{\it negative}) if all
exponents of $A$ in $W$ are positive (negative).

We are interested in the matrices that result when the two letters
are (independent) positive definite (complex Hermitian) $n$-by-$n$
matrices (PD, for short).  For convenience, the letters $A,B$ will
also represent the substituted PD matrices (the context will make
the distinction clear).  To make sure that $W$ is well-defined
after substitution, we take primary PD powers (see \cite[p.
433]{HJ2} and \cite[p. 413]{HJ2}). I.e. given $p \in \mathbb
R\backslash\{0\}$, a unitary matrix $U$, and a nonnegative
diagonal matrix $D$, we have ${(UDU^{*})}^{p} = UD^{p}U^{*}$.

Our primary interest will be matrix equations involving
$A$-positive symmetric g-words.

\begin{definition}
A {\it symmetric word equation} is an equation, $S(A,B) = P$, in
which $S(A,B)$ is an $A$-positive symmetric g-word. If $B$ and $P$
are given positive definite matrices, any positive definite matrix
$A$ for which the equation holds is called a {\it solution} to the
symmetric word equation.
\end{definition}

A symmetric word equation will be called {\it solvable} if there
exists a solution for every pair of positive definite $n$-by-$n$
$B$,$P$. Moreover, if each such pair gives rise to a unique
solution, the equation will be called {\it uniquely solvable}. It
is clear that the (unique) solvability of $A$-positive word
equations implies the same as for $A$-negative equations (by
replacing $A$ with $A^{-1}$), and hence, no generality is lost in
Definition 1.1. As we shall soon see, the restriction of the
definition above to $A$-positive words is important.

We first encountered symmetric word equations while studying a
trace conjecture \cite{JH} in the case of words with positive
integral powers.  The conjecture is

\begin{conjecture}
A word has positive trace for every pair of real positive definite
matrices if and only if the word is symmetric or a product
(juxtaposition) of 2 symmetric words.
\end{conjecture}

It turns out that for each solvable symmetric word equation, we
can identify an infinite class of words that admit real PD
matrices $A$ and $B$ giving those words a negative trace.  Our
investigation of these equations, although useful for our methods
in \cite{JH}, show that they have a rich structure of their own,
some of which we explore here.  (See Section 7 for a
generalization of the notion of "symmetric word equation" defined
above).  A natural question to ask is if it is necessary to focus
attention on symmetric g-words.  We remark that it can be shown
that symmetric g-words are the only g-words that are positive
definite for all positive definite matrices $A,B$.  In light of
Definition 1.1, our restriction seems appropriate.

\section{Symmetric words}

Recall that two $n$-by-$n$ matrices $X$ and $Y$ are said to be
{\it congruent} if there is an invertible $n$-by-$n$ matrix $Z$
such that $Y = Z^{*}XZ$; and that congruence on Hermitian matrices
preserves inertia (the ordered triple consisting of the number of
positive, negative, and zero eigenvalues) and, thus, positive
definiteness \cite[p. 223]{HJ1}. A symmetric word of class $k$ in
two positive definite matrices is congruent to one of class $k-1$,
iteration of which implies congruence to the "center," class $0$,
positive definite matrix. We conclude that

\begin{lemma}
A symmetric g-word in two positive definite matrices is positive
definite.
\end{lemma}

A natural question to ask is if the map from the set of PD
matrices to itself given by $h: A \mapsto S(A,B)$ is surjective.
Our main result answers this in the affirmative.

\begin{theorem}
Every symmetric word equation is solvable.
\end{theorem}

We remark that the statement is not generally valid if the
definition of "symmetric word equation" is widened to allow mixed
sign powers of $A$, even in the case $n > 1$.  Let $I$ denote the
$n$-by-$n$ identity matrix and let $B \neq I$ be an $n$-by-$n$ PD
matrix. Then, the equation, $A^{-1}BA^{2}BA^{-1} = I$, has no PD
solution. For if there were a PD solution $A$, then $ABA^{-1} =
U$, for some unitary $U$.  Hence, $B = A^{-1}UA$ has eigenvalues
on the unit disc.  The only positive definite matrix for which
this is so is $B = I$.

We believe that $h$ is also injective, but this is proven only in
special cases.  This is

\begin{conjecture}
Every symmetric word equation is uniquely solvable.
\end{conjecture}

\section{Some specific equations}

In what follows, we shall say that two symmetric word equations
are ({\it uniquely}) {\it equivalent} if for each, its (unique)
solvability implies the (unique) solvability of the other.  For
instance, the equations $ABA^{2}BA = P$ and $ABA = P$ are uniquely
equivalent because PD matrices have unique PD square roots.  As a
more subtle example, the equations $A^{2}BABA^{2} = P$ and
$ABA^{3}BA = P$ are seen to be uniquely equivalent by setting $A =
B^{-1/2}XB^{-1/2}$ in $A^{2}BABA^{2} = P$. Additionally, the same
substitution gives us that for any integer $k > 2$,
$A^{2}(BA)^{k-2}BA^{2} = P$ and $ABA^{k}BA = P$ are uniquely
equivalent.

It is an easy exercise that the facts below follow from our
definition of unique equivalence:
\begin{enumerate}
\item For $r > 0$ and $s \neq 0$, the equations $S(A^{r},B^{s}) =
P$ and $S(A,B) = P$ are uniquely equivalent.

\item For each positive integer $k$, $S(A,B)^{k} = P$ and $S(A,B) =
P$ are uniquely equivalent.

\item For $s$ real, the equations $B^{s}S(A,B)B^{s} = P$ and
$S(A,B) = P$ are uniquely equivalent.
\end{enumerate}
The simplest examples of symmetric word equations are $\{A^{q} =
P,\ q\in \mathbb R\backslash \{ 0\}  \}$.  These satisfy
Conjecture 2.3 by uniqueness of PD $q^{th}$ roots. More
interesting is the first non-trivial equation, $ABA = P$. This
equation has arisen in other contexts \cite{Ando}.  We indicate
some aspects of this equation useful to us.

\begin{theorem}
The equation, $ABA = P$, has a unique solution for each pair of PD
matrices $B,P$.  Moreover, the unique PD matrix $A$ is given by
\[A = B^{-1/2}(B^{1/2}PB^{1/2})^{1/2}B^{-1/2}.\]
\end{theorem}

\begin{proof}
Assume that $P$ and $B$ are given positive definite matrices and
$A$ is a PD solution to $ABA = P$.  Set $X = B^{1/2}AB^{1/2}$,
which is PD by Lemma 2.1, so that $A = B^{-1/2}XB^{-1/2}$. Then,
\[P = B^{-1/2}XB^{-1/2}BB^{-1/2}XB^{-1/2} =
B^{-1/2}X^{2}B^{-1/2}.\] Therefore, $X^{2} = B^{1/2}PB^{1/2}$,
from which it follows that $X$ is uniquely determined as
$(B^{1/2}PB^{1/2})^{1/2}$. Hence, $A$ must be
$B^{-1/2}(B^{1/2}PB^{1/2})^{1/2}B^{-1/2}$. Finally, substituting
this positive definite $A$ (by Lemma 2.1 again) into the original
equation does verify that it is a solution.
\end{proof}

As in \cite{Ando}, given two PD matrices $C$ and $D$ we will
denote $C\#D$ as the PD matrix, \[C\#D =
C^{1/2}(C^{-1/2}DC^{-1/2})^{1/2}C^{1/2},\] the so-called {\it
geometric mean} of $C$ and $D$.  Notice that from Theorem 3.1,
$ABA = P$ has the unique solution $A = (B^{-1})\#P$.  Writing this
equation (by inverting) as $B^{-1} = AP^{-1}A$, we also have that
$A = P\#(B^{-1})$, and from this comes the not so obvious fact
that $C\#D = D\#C$.

\section{Fixed points and a sequence lemma}

For more complex symmetric word equations, it is not clear that
there should be an explicit formula for a solution, as in Theorem
3.1. We do not know one, for example, for $A^{2}BABA^{2} = P$.
Since our proof of the solvability of this equation (and all
others) will use fixed-point theory, we record a useful theorem of
Brouwer \cite{Zeid}.

\begin{theorem}[Brouwer's fixed point theorem]
If $M$ is a compact, convex subset of a finite dimensional Banach
space and if $f : M \to M$ is a continuous function, then there is
a fixed point, $p$, for $f$ in $M$.
\end{theorem}

We will be using the spectral matrix norm throughout (see \cite[p.
295]{HJ1}).  This norm is useful because for positive semidefinite
$A$, it is just the largest eigenvalue of $A$. Brouwer's result is
an important ingredient in the proof of Theorem 2.2. Before
proving Theorem 2.2, we record the following.

\begin{lemma}
Suppose $\{A_{k}\}$ is a convergent sequence of positive definite
matrices.  Then, there is a subsequence $\left\{ {A_{k_j } }
\right\}_{j = 1}^\infty$ such that \[ \mathop {\lim }\limits_{j
\to \infty } \;\;\frac{{(A_{k_j } + I/k_j )^{ - 1} }}{{\left\|
{\left( {A_{k_j }  + I/k_j } \right)^{ - 1} } \right\|}}\] exists.
\end{lemma}

\begin{proof}
Set $A_{k} = U_{k}D_{k}U_{k}^{*}$ for unitary $U_{k}$ and $D_{k}
=$ diag$(\lambda_{1k}, \lambda_{2k} , \ldots, \lambda_{nk})$ in
which $0 < \lambda_{1k} \leq \lambda_{2k} \leq \ldots \leq
\lambda_{nk}$ are the eigenvalues of $A_{k}$. For each $m \in
\{1,\ldots,n\}$, define new sequences, \[h_{m}(k) = \frac{{\lambda
_{1k} + 1/k}}{{ \lambda _{mk} + 1/k}}.\] Notice that $h_{1}(k) =
1$ for all $k$.  We now show by induction that there exists a
subsequence, $\left\{ {k_j } \right\}_{j = 1}^\infty$, such that
each of $\left\{ {h_{m}(k_j) } \right\}_{j = 1}^\infty$ converges.
To simplify matters later, we first assume that $\{U_{k}\}$
converges to a unitary $U$ by passing to a subsequence (the set of
unitary matrices is compact).

Let $M \in \{1,\ldots,n\}$ be such that there exists a subsequence
$\left\{ {k_j } \right\}_{j = 1}^\infty$ making each of $\left\{
{h_{1}(k_j) } \right\}_{j = 1}^\infty$,$\ldots$,$\left\{
{h_{M}(k_j) } \right\}_{j = 1}^\infty$ converge.  Such an $M$
clearly exists (e.g. $M = 1$).  If $M = n$, then there is nothing
to prove.  Otherwise, examine the inequality, \[0 < h_{M+1}(k_j) =
\frac{{\lambda _{1k_j} + 1/k_j}}{{ \lambda _{(M+1)k_j} + 1/k_j}}
\leq 1,\] which holds for all $j$.  Since this is a bounded
sequence, choose a subsequence of the sequence, $\left\{ {k_j }
\right\}_{j = 1}^\infty$, making this ratio converge to some
nonnegative number. This will not alter the convergence of the
first $M$ sequences. This completes the induction.

We conclude that there exists a subsequence, $\left\{ {k_j }
\right\}_{j = 1}^\infty$, such that \[\mathop {\lim }\limits_{j
\to \infty } {\kern 1pt} {\kern 1pt} {\kern 1pt} \frac{{ \lambda
_{1k_j } + 1/k_j }}{{ \lambda _{mk_j } + 1/k_j }} = h_m \ge 0.\]
We now claim that $\mathop {\lim }\limits_{j \to \infty }
\frac{{(A_{k_j }  + I/k_j )^{ - 1} }}{{\left\| {\left( {A_{k_j } +
I/k_j } \right)^{ - 1} } \right\|}}$ exists.  But the matrix in
question is just \[U_{k_j } \left( {\begin{matrix}
   {\frac{\lambda _{1k_j }  + 1/k_j }
{\lambda _{1k_j }  + 1/k_j }} &  \ldots  & 0  \\
    \vdots  &  \ddots  &  \vdots   \\
   0 &  \cdots  & {\frac{\lambda _{1k_j }  + 1/k_j }
{\lambda _{nk_j }  + 1/k_j }}  \\
 \end{matrix} } \right)U_{k_j } ^* \] which by construction has the limit, \[U\left(
 {\begin{matrix}
   1 &  \ldots  & 0  \\
    \vdots  &  \ddots  &  \vdots   \\
   0 &  \cdots  & {h_n }  \\
 \end{matrix}} \right)U^*.\] This proves the lemma.
\end{proof}

\section{All symmetric word equations are solvable}

We may now prove the main result.

\begin{proof}[Proof of Theorem 2.2]
Notice (using (1), (3) above) that it is equivalent to study
symmetric word equations of the form $AS(A,B)A = P$, where
$S(A,B)$ is a symmetric g-word beginning and ending with a power
of $A$.  For instance, the equation $A^{1/2}BABA^{1/2} = P$ is
(uniquely) equivalent to $A(ABA^{4}BA)A = P$.  The convex, compact
set to which we apply Brouwer's theorem is \[M = \{A \ | \  A
\text{ is positive semidefinite and } \|A\| \leq 1\}\]

Define $S_{A}$ to denote the sum of all the powers of $A$ in
$S(A,B)$. Also, let $S_{-B}$ $(S_{B})$ denote the sum of all
negative (positive) powers of $B$ in $S(A,B)$.  Let $k$ be a
positive integer. Define a function $f_{k}$ on the set $M$ as
follows:
\begin{equation*}
\begin{split}
f_k (X) &= \frac{{P\# \left( {S\left( {X + I/k,B} \right)^{ - 1} }
\right)}}{{g_k (X)}}\\ &= \frac{{P^{1/2} \left( {P^{ - 1/2} S(X +
I/k,B)^{ - 1} P^{ - 1/2} } \right)^{1/2} P^{1/2} }}{{g_k (X)}}
\end{split}
\end{equation*} where $g_{k}: M \to \mathbb R^{+}$ is defined by, \[ g_{k} (X) =
\left\| P \right\|\left\| {P^{ - 1} } \right\|^{1/2} \left\| B
\right\|^{S_{ - B} /2} \left\| {B^{ - 1} } \right\|^{S_B /2}
\left\| {\left( {X + I/k} \right)^{ - 1} } \right\|^{S_A /2}.\]

From the properties of the spectral matrix norm and the fact that
$S(C,D)^{-1} = S(C^{-1},D^{-1})$, it follows that $f_{k}$ is
bounded by 1. I.e., \[ \left\| {f_k (X){\kern 1pt} } \right\| \le
\frac{{\left\| P \right\|\left\| {P^{ - 1} } \right\|^{1/2}
\left\| B \right\|^{S_{ - B} /2} \left\| {B^{ - 1} } \right\|^{S_B
/2} \left\| {\left( {X + I/k} \right)^{ - 1} } \right\|^{S_A /2}
}}{{g_k (X)}} = 1.\]

We should note that in the 1-by-1 case, we simply have $f_{k}(X) =
1$, so that the unique fixed point for $f_{k}(X)$ is given by $X =
1$. More generally, from the properties of the geometric mean and
since $k > 0$, it is also clear that $f_{k}(X)$ is positive
definite.  Hence, $f_{k}(X)\in M$. From the discussion of primary
matrix functions \cite[p. 433]{HJ2} for normal matrices, it is
seen that for $q \neq 0$, $X^{q}$ is a continuous function on $M$.
Therefore, $f_{k}$ is also continuous on $M$ (it is made up of
compositions and products of continuous functions).

Now, apply Brouwer's fixed point theorem to give us $A_{k} =
f_{k}(A_{k})$. Because $f_{k}$ is always positive definite, we
must have that $A_{k}$ is nonsingular.  Hence, from the properties
of the geometric mean (Theorem 3.1), we have
\begin{equation}
g_{k}(A_{k})^{2}A_{k} S(A_{k} + I/k,B)A_{k} = P.
\end{equation}

Because $\left\{ {A_{k}} \right\}_{k = 1}^\infty$ is an infinite,
bounded sequence of PD matrices, there is a subsequence that
converges to a positive semidefinite matrix $T$.  We will
therefore assume that $\{A_{k}\}$ converges.  If $T$ is actually
positive definite, then we have \[P = \mathop {\lim }\limits_{k
\to \infty } \,\,g_k (A_k )^2 A_k {\kern 1pt} S\left( {A_k  +
I/k,B} \right)A_k  = TS\left( {T,B} \right)T\,g_\infty  (T)^2\] in
which $g_\infty  (T) = \left\| P \right\|\left\| {P^{ - 1} }
\right\|^{1/2} \left\| B \right\|^{S_{ - B} /2} \left\| {B^{ - 1}
} \right\|^{S_B /2} \left\| {T^{ - 1} } \right\|^{S_A /2} $ .
Then, it is easily seen that \[A = T \cdot g_\infty  (T)^{2/\left(
{2 + S_A } \right)}\] is our desired PD solution.

Since (5.1) has no limit interpretation if $T$ is singular, we now
show that $T$ is necessarily positive definite. For each $k$, form
the decomposition, $A_{k} = U_{k}D_{k}U_{k}^{*}$, with unitary
$U_{k}$ and $D_{k} =$ diag$(\lambda_{1k}, \lambda_{2k} , \ldots,
\lambda_{nk})$, where $0 < \lambda_{1k} \leq \lambda_{2k} \leq
\ldots\leq \lambda_{nk}$ are the eigenvalues of $A_{k}$. Now,
suppose that $\{A_{k}\}$ converges to a singular matrix $T$ with
$0 \leq \lambda_{1} \leq \lambda_{2} \leq \ldots\leq \lambda_{n}$
being the eigenvalues of $T$. Since $T$ is assumed to be singular,
let $m \geq 1$ be such that
$\{\lambda_{1},\lambda_{2},\ldots,\lambda_{m}\}$ are the zero
eigenvalues of $T$. Set $c = \left\| P \right\|^2 \left\| {P^{ -
1} } \right\|\left\| B \right\|^{S_{ - B} } \left\| {B^{ - 1} }
\right\|^{S_B }$, and examine the equality (following from
inverting equation (5.1)) \[ \left\| {\left( {A_k  + I/k}
\right)^{ - 1} } \right\|^{ - S_A } S((A_k  + I/k)^{ - 1} ,B^{ -
1} ) = cA_k P^{ - 1} A_k \] which is just
\begin{equation}
S\left( {\frac{{(A_k  + I/k)^{ - 1} }}{{\left\| {\left( {A_k  +
I/k} \right)^{ - 1} } \right\|}},B^{ - 1} } \right) = cA_k P^{ -
1} A_k.
\end{equation}

From Lemma 4.2, there is a subsequence $\left\{ {A_{k_j } }
\right\}_{j = 1}^\infty$ such that $\left\{ {U_{k_j } }
\right\}_{j = 1}^\infty$ converges to some unitary $U$ and such
that $\mathop {\lim }\limits_{j \to \infty } \frac{{(A_{k_j }  +
I/k_j )^{ - r} }}{{\left\| {\left( {A_{k_j }  + I/k_j } \right)^{
- 1} } \right\|^r }}$ exists for all $r > 0$.  Moreover, this
limit is equal to
\begin{equation}
U\left( {\begin{matrix}
   1 &  \ldots  & 0  \\
    \vdots  &  \ddots  &  \vdots   \\
   0 &  \cdots  & {h_n ^r }  \\
 \end{matrix} } \right)U^*
\end{equation} for some $h_{2},h_{3},\ldots,h_{n} \geq 0$. If
$h_{2},h_{3},\ldots,h_{n}$ are all nonzero, then the limit as
$j\rightarrow \infty$ of \[S\left( {\frac{{(A_{k_j }  + I/k_j)^{ -
1} }}{{\left\| {\left( {A_{k_j } + I/k_j} \right)^{ - 1} }
\right\|}},B^{ - 1} } \right)\] is invertible, while the limit of
$cA_{k_j } P^{ - 1} A_{k_j }$ is singular. Whence, $h_{t} = 0$ for
some $t\in\{2,3,\ldots,n\}$. Moreover, if $s > t$, then $h_{s} =
0$ as well; this coming from the fact that \[h_{t}(k_{j}) =
\frac{{ \lambda _{1k_j } + 1/k_j }}{{ \lambda _{t{\kern 1pt} k_j }
+ 1/k_j }} \ge \frac{{ \lambda _{1k_j } + 1/k_j }}{{ \lambda
_{s{\kern 1pt} k_j } + 1/k_j }} = h_{s}(k_{j})\] by our ordering
of the eigenvalues of $A_{k}$.  So assume that $t$ is the largest
element of $\{1,2,\ldots,n\}$ such that $h_{t} \neq 0$. Then, our
limit (5.3) looks like \[L = U\left( {\begin{matrix}
   {E_h ^r } & 0  \\
   0 & 0  \\
 \end{matrix} } \right)U^*  = UE^rU^*\] in which
$E_{h}$ = diag$(1,h_{2},\ldots,h_{t})$ is positive and $E$ =
diag$(1,h_{2},\ldots,h_{t},0,\ldots,0)$. Now, set $\widetilde B =
U^* B^{ - 1} U$ .  Then, the left hand side of (5.2) (within the
subsequence above) converges to
\begin{equation}
S\left( UEU^{*},B^{-1} \right) = UE^{p_1 } \widetilde B^{q_1 }
E^{p_2 } \widetilde B^{q_2 } \cdots \widetilde B^{q_2 } E^{p_2 }
\widetilde B^{q_1 } E^{p_1 } U^*.
\end{equation}

We claim (5.4) has the form \[U\left( {\begin{matrix}
   H & 0  \\
   0 & 0  \\
 \end{matrix} } \right)U^*\] for some PD $t$-by-$t$ matrix $H$.
Indeed, the center matrix in (5.4) is of the form $E^p $ or $E^p
\widetilde B^q E^p$, each being the direct sum of a $t$-by-$t$ PD
matrix and an ($n-t$)-by-($n-t$) zero matrix.  We now induct on
the form of (5.4). Assume that $E^{p_d } \widetilde B^{q_d }
\cdots \widetilde B^{q_d } E^{p_d }$ is a direct sum of a
$t$-by-$t$ PD matrix $H_{d}$ with a zero matrix.  Form the
partition,
\begin{equation}
\widetilde B^{q_{d - 1} }  = \left( {\begin{matrix}
   {B_{11} } & {B_{12} }  \\
   {B_{21} } & {B_{22} }  \\
 \end{matrix} } \right)
\end{equation} in which $B_{11}$ is a $t$-by-$t$ PD matrix, $B_{22}$ is an
($n-t$)-by-($n-t$) PD matrix and $B_{12} = B_{21}^{*}$ (see
\cite[p. 472]{HJ1}). Then,
\begin{equation*}
\begin{split}
E^{p_{d - 1} } \widetilde B^{q_{d - 1} } \left( {\begin{matrix}
   {H_d } & 0  \\
   0 & 0  \\
 \end{matrix} } \right)\widetilde B^{q_{d - 1} } E^{p_{d - 1} } &=
E^{p_{d - 1} } \left( {\begin{matrix}
   {B_{11} } & {B_{12} }  \\
   {B_{21} } & {B_{22} }  \\
 \end{matrix} } \right)\left( {\begin{matrix}
   {H_d } & 0  \\
   0 & 0  \\
 \end{matrix} } \right)\left( {\begin{matrix}
   {B_{11} } & {B_{12} }  \\
   {B_{21} } & {B_{22} }  \\
 \end{matrix} } \right)E^{p_{d - 1} }\\
&=E^{p_{d - 1} } \left( {\begin{matrix}
   {B_{11} H_d B_{11} } & {B_{11} H_d B_{12} }  \\
   {B_{21} H_d B_{11} } & {B_{21} H_d B_{12} }  \\
 \end{matrix} } \right)E^{p_{d - 1} }\\
&= \left( {\begin{matrix}
   {E_h ^{p_{d - 1} } B_{11} H_d B_{11} E_h ^{p_{d - 1} } } & 0  \\
   0 & 0  \\
 \end{matrix} } \right)
\end{split}
\end{equation*} is of the desired form, completing the induction.

Examine now the right hand side of (5.2), $cA_{k_j } P^{ - 1}
A_{k_j }$, which converges to $cTP^{-1} T$.  Since the left hand
side limit by above has rank $t > 0$, $T$ cannot have all its
eigenvalues equal to 0.  Whence, we can write this limit as
$cUDU^* P^{ - 1} UDU^*$, in which $D$ = $\left( {\begin{matrix}
   0 & 0  \\
   0 & \Lambda   \\
 \end{matrix} } \right)$ and $\Lambda$ is the positive diagonal matrix,
diag$(\lambda_{m+1},\ldots,\lambda_{n})$. Setting these two
expressions equal gives us
\[c\left( {\begin{matrix}
   0 & 0  \\
   0 & \Lambda   \\
 \end{matrix} } \right)U^* P^{ - 1} U\left( {\begin{matrix}
   0 & 0  \\
   0 & \Lambda   \\
 \end{matrix} } \right) = \left( {\begin{matrix}
   H & 0  \\
   0 & 0  \\
 \end{matrix} } \right).\]

Of course, this equality is impossible. We simply compare the two
$(n,n)$ entries of the left and right hand sides to arrive at a
contradiction (one is positive, the other zero). Hence, $T$ is
invertible, completing the proof.
\end{proof}

We remark that if $P$ and $B$ are chosen to be real, then the
proof above shows that the solution $A$ may be chosen to be real
as well.

\section{Approximate solutions}

We now make some remarks about finding approximate solutions to
symmetric word equations.  As a first approach one would hope that
an iteration of the function defined in the proof of Theorem 2.2
would give rise to approximate solutions.  Unfortunately,
experimentation shows this not to be the case.  We detail another
method that has been quite effective in practice at finding
solutions (in the case of positive integral powers) and in
verifying our conjecture of unique solvability.

Given a symmetric word $S(A,B)$ with positive integral exponents,
start with an initial PD matrix, $A_{0}$, (usually $I$), and
expand the expression, $S(A_{0}+D,B)$.  Consider the formal sum,
$S'(A_{0}+D,B)$, of the terms in this expansion with at most a
single $D$. Now, solve the linear system \[S'(A_{0}+D,B) = P\] for
the matrix $D$ and set $A_{1}$ $ \leftarrow$ $A_{0} + D$.
Repeating this process gives our algorithm.  As a simple example,
the repeated equations for $ABA = P$ are given by
\[A_{i-1}BD + DBA_{i-1} = P - A_{i-1}BA_{i-1}\] \[A_{i} = A_{i-1}
+ D.\]  Curiously enough, there seems to be no guarantee in
general that these $A_{i}$ will be positive definite (or even
Hermitian), nor is it clear that the linear system for $D$ above
will always have a solution.  Nonetheless, experimentation has
shown that these iterations always converge to the same PD
solution regardless of initial starting point.

\section{A generalization}

We close by noting a generalization of Theorem 2.2 to a larger
class of matrix equations.  If $C_1,\ldots,C_m$ is any list of $m$
invertible $n$-by-$n$ matrices and $W(A; C_1,\ldots,C_m) =
A^{p_{1}}C_{1}A^{p_{2}}C_{2} \cdots
C_{2}^{*}A^{p_{2}}C_{1}^{*}A^{p_{1}}$ is such that $W(A;
C_1,\ldots,C_m)^*$ $= W(A; C_1,\ldots,C_m)$, we call $W$ a {\it
generalized symmetric word}.  It is again an elementary exercise
in congruence that, if $A$ is PD and $W$ is a generalized
symmetric word, then $W(A; C_1,\ldots,C_m)$ is positive definite.
For a given PD matrix $P$ and invertible $C_1,\ldots,C_m$ we call
\[W(A; C_1,\ldots,C_m) = P\] with each $p_i > 0$, a {\it
generalized symmetric word equation} in the PD variable $A$.
Though we are not motivated by any particular application,
generalized symmetric word equations are natural to consider on a
theoretical level.  Our proof of Theorem 2.2 goes over directly to
generalized symmetric word equations, except that a further
technical condition on the $C_i$'s is needed.  That condition is
that the leading principal minors of any unitary similarity of
$C_i^{-1}$ should be nonzero. Of course, this means that any
principal minor of any unitary similarity of $C_i^{-1}$ should be
nonzero for $i=1,\ldots,m$. Using elementary facts about the field
of values of an $n$-by-$n$ matrix $C$: $$F(C) = \{x^*Cx \ | \ x^*x
= 1, \ x \in \mathbb C^n\},$$ see chapter 1 of \cite{HJ2}, the
latter condition may easily be seen to be equivalent to the
statement that $0 \notin F(C)$.  For purposes of this work, we
call such a matrix {\it completely invertible}. Of course, any
power of a positive definite matrix $B$ is completely invertible,
and this is what is essential in the proof of Theorem 2.2.  For
generalized symmetric word equations, the analysis of (5.4)
remains valid because any partition as in (5.5) will have an
invertible leading principle submatrix. We conclude

\begin{theorem}
For completely invertible $C_1,\ldots,C_m$, the generalized
symmetric word equation \[W(A; C_1,\ldots,C_m) = P\] is solvable
for any positive definite $P$.
\end{theorem}

We have no example showing that complete invertibility cannot be
replaced by invertibility.

\section{Acknowledgement}

The authors would like to thank Scott Armstrong for a careful
reading of a preliminary version of this manuscript.

\bibliographystyle{amsplain}

\end{document}